\documentclass[reqno,centertags,12pt]{amsart}
\title[Inhomogeneous monostable transition fronts]{Transition fronts for inhomogeneous monostable reaction-diffusion equations via linearization at zero}
\author{Tianyu Tao, Beite Zhu, and Andrej Zlato\v s}

\address{\noindent University of
Wisconsin \\ Madison, WI 53706, USA \\ Email: \rm
ttao@wisc.edu}

\address{\noindent University of
California \\ Berkeley, CA 94720, USA \\ Email: \rm
jupiter\_ju@berkeley.edu}

\address{\noindent Department of Mathematics \\ University of
Wisconsin \\ Madison, WI 53706, USA \newline Email: \rm
zlatos@math.wisc.edu}

\usepackage{amsmath,amsthm,amscd,amssymb,latexsym,verbatim} 
\usepackage{setspace}
\usepackage{graphicx,epsf,cite}

-\marginparwidth \oddsidemargin -4mm \evensidemargin \oddsidemargin

\newtheorem{theorem}{Theorem}[section]

\newtheorem*{theorem*}{Theorem}

\newtheorem{lemma}[theorem]{Lemma}

\theoremstyle{definition}

\theoremstyle{remark}

\newcounter{smalllist}

\DeclareMathOperator{\diam}{diam}

\allowdisplaybreaks
\numberwithin{equation}{section}

\newcommand{\lb}{\label}

\newcommand{\supp}{\text{\rm{supp}}}

\newcommand{\beq}{\begin{equation}}
\newcommand{\eeq}{\end{equation}}

\newcommand{\bal}{\begin{align}}
\newcommand{\eal}{\end{align}}
\newcommand{\bals}{\begin{align*}}
\newcommand{\eals}{\end{align*}}

\newcommand{\eps}{\varepsilon}

\newcommand{\tht}{\theta}

\newcommand{\til}{\tilde}

\newcommand{\bbR}{{\mathbb{R}}}

\usepackage{framed}

\begin{document}

\begin{abstract}
We prove existence of transition fronts for a large class of reaction-diffusion equations in one dimension, with inhomogeneous monostable reactions.  We construct these as perturbations of corresponding front-like solutions to the linearization of the PDE at $u=0$.  While a close relationship of the solutions to the two PDEs has been well known and exploited for KPP reactions (and our method is an extension of such ideas from \cite{ZlaInhomog}), to the best of our knowledge this is the first time such an approach has been used in the construction and study of  fronts for non-KPP monostable reactions. 
\end{abstract}
\maketitle

\section{Introduction}

We study transition fronts for the one-dimensional reaction-diffusion equation
\begin{equation} \label{1.1}
u_t=u_{xx}+f(x,u),
\end{equation}
with an inhomogeneous non-negative reaction $f\ge 0$ satisfying $f(x,0)=f(x,1)=0$, and with $u\in[0,1]$.  Such PDEs model a host of natural processes such as combustion, chemical reactions, population dynamics and others, with  $u$ representing (normalized) temperature, concentration of a reactant, or population density.  

Both $u\equiv 0$ and $u\equiv 1$ are equilibrium solutions of \eqref{1.1} and one is interested in the study of propagation of reaction in space, that is,  invasion of the state $u=0$ by the state $u=1$.  An important class of solutions modeling the propagation of reaction are transition fronts.   A {\it (right-moving) transition front} is any entire solution $u:\bbR^2\to[0,1]$ of \eqref{1.1}  which satisfies
\begin{equation} \label{1.6}
\lim_{x\to -\infty} u(t,x) =1 \qquad \text{ and } \qquad \lim_{x\to +\infty} u(t,x) =0
\end{equation} 
for each $ t\in \mathbb{R}$.  In addition, we also require that for any $\eps > 0$ there exists $L_{\eps} < \infty$ such that
\begin{equation} \label{1.7}
\sup_{t \in \bbR}\diam\{x \in \bbR \,|\, \eps \le u(t,x)\le 1-\eps\} \le L_{\eps}. 
\end{equation}
The definition of a left-moving transition front is similar, with the limits in \eqref{1.6} exchanged.  We will only study right-moving fronts here because the treatment of both cases is identical, up to a reflection in $x$.  We note that the above definition is from \cite{BH3,Shen,Matano}.

We will consider  here the case of {\it monostable reactions}, for which $u\equiv1$ is an asymptotically stable solution while $u\equiv 0$ is unstable.  We assume that $f$ is Lipschitz, 
\beq\lb{1.30}
f(x,0)=f(x,1) = 0 \quad \text{ for } x\in\bbR,
\eeq 
\beq\lb{1.39}
a(x) := f_u(x,0) > 0
\eeq
exists, and
\begin{equation} \label{1.2}
a(x)g_0(u) \le f(x,u) \le a(x)g_1(u) \quad \quad \text{ for } (x,u) \in \mathbb{R}\times [0,1],
\end{equation}
where $g_0,g_1 \in C^1([0,1])$  satisfy
\begin{alignat}{4}
\label{1.3.1}&g_0(0)=g_0(1)=0,\quad\quad &g_0'(0)=1,\quad \quad &g_0(u)>0 \quad \text{ and } \quad g_0'(u)\le 1 \quad \text{ for } u \in (0,1),&\\ 
\label{1.3.2}&g_1(0)=0,\quad &g_1'(0)=1,\quad\quad &g_1'(u)\ge 1 \quad \text{ for } u\in[0,1], &
\end{alignat}
\begin{equation} \label{1.4}
\int_0^1 \frac{g_1(u)-g_0(u)}{u^2}du <\infty.
\end{equation}
Finally, we let
\beq \lb{1.34}
a_{-}:=\inf_{x \in \mathbb{R}}a(x) \le \sup_{x \in \mathbb{R}}a(x) =: a_+.
\eeq

When the reaction $f(x,u)=f(u)\ge 0$ is homogeneous, a special case of transition fronts are {\it traveling fronts}.  These are  of the form $u(t,x)=U(x-ct)$, with some front speed $c$ and front profile $U$ such that $\lim_{s\to-\infty} U(s)=1$ and $\lim_{s\to\infty} U(s)=0$, and  their study goes back to the seminal works of Kolmogorov, Petrovskii, and Piskunov \cite{KPP}, and Fisher \cite{Fisher}.  They considered {\it KPP reactions}, a special case of monostable reactions with $g_1(u)=u$, and found that for each $c\ge c_0:=2\sqrt {f'(0)}$ there is a unique traveling front $u(t,x)=U_c(x-ct)$.  A simple phase-plane analysis shows that this turns out to be the case for general homogeneous monostable reactions, although with a different $c_0\ge 2\sqrt {f'(0)}$.  In contrast, {\it ignition reactions}, satisfying $f(u)=0$ for $u\in[0,\tht]\cup\{1\}$ and $f(u)>0$ for $u\in(\tht,1)$ (for some {\it ignition temperature} $\tht\in(0,1)$), give rise to a single speed $c_0>0$ and a single traveling front. 

Despite many developments for homogeneous and space-periodic reactions in the almost eight decades since \cite{KPP,Fisher} (see the reviews \cite{Berrev,Xin2} and references therein),  transition fronts in spatially non-periodic media have only been studied relatively recently. The first existence result, for small perturbations of homogeneous {\it bistable reactions} (the latter are such that $f(u)<0$ for $u\in(0,\tht)$ and $f(u)>0$ for $u\in(\tht,1)$), was obtained by Vakulenko and Volpert \cite{VakVol}.  Existence results without a hypothesis of closeness to a homogeneous reaction, for ignition reactions of the form $f(x,u)=a(x)g(u)$ for some homogeneous ignition reaction $g$, were proved by Mellet,  Roquejoffre, and Sire \cite{MRS}, and by Nolen and Ryzhik \cite{NolRyz} (see also \cite{MNRR} for uniqueness and stability results for these reactions).  Existence of fronts for general inhomogeneous ignition reactions as well as for some monostable reactions which are in some sense not too far from ignition ones was proved by Zlato\v s \cite{ZlaGenfronts} (uniqueness and stability for the ignition case was also obtained).  All these results are based on recovering a front as a locally uniform limit, along a subsequence, of solutions $u_n$ of the Cauchy problem with initial data $u_n(\tau_n,x)\approx\chi_{(-\infty,-n)}(x)$, where $\tau_n\to-\infty$ are such that $u_n(0,0)=\tfrac 12$.  Existence of a limit $u$ on $\bbR^2$ is guaranteed by  parabolic regularity, and the challenge is to show that $u$ is a transition front.  We note that even in the monostable case in \cite{ZlaInhomog}, when one expects multiple transition fronts,  existence of only a single transition front was obtained. 

A very different approach has been used by Nolen, Roquejoffre, Ryzhik, and Zlato\v s \cite{NRRZ}, and by Zlato\v s \cite{ZlaInhomog} to prove existence of multiple transition fronts for inhomogeneous KPP reactions.  It is well known that when $f$ is KPP, then there is a close relationship between the solutions of \eqref{1.1} and those of its linearization 
\begin{equation}  \label{1.8}
v_t = v_{xx} +a(x)v
\end{equation}
at $u=0$.  The reason for this is that all KPP fronts are {\it pulled}, with the front speeds determined by the reaction at $u=0$, which is due to the {\it reaction strength} $\tfrac{f(x,u)}u$ being largest at $u=0$ for any fixed $x\in\bbR$.  This is in stark contrast with ignition fronts, which are always {\it pushed} because they are ``driven'' by the reaction at intermediate values of $u$.  

One can therefore consider the simpler {\it front-like solutions} of \eqref{1.8}, which
are of the form 
\beq\lb{1.31}
v_{\lambda}(t,x) = e^{\lambda t}\phi_{\lambda}(x).
\eeq
Here $\phi_{\lambda}>0$ is a generalized eigenfunction of the operator $\mathcal{L}:=\partial_{xx}+a(x)$, satisfying 
\begin{equation} \label{1.9}
\phi''_{\lambda} + a(x) \phi_{\lambda} = \lambda \phi_{\lambda}
\end{equation}
on $\bbR$, which  exponentially grows to $\infty$ as $x\to-\infty$ and exponentially  decays to 0 as $x\to\infty$.  If we let $\lambda_0:=\sup \sigma(\mathcal L)\in[a_-,a_+]$ be the supremum of the spectrum of $\mathcal L$, then it is well known that such $\phi_\lambda$ exists precisely when $\lambda>\lambda_0$, and is unique if we also require $\phi_\lambda(0)=1$.

For KPP reactions one can try to use these solutions to find transition fronts for \eqref{1.1} with
\beq\lb{1.41}
\lim_{x\to\infty} \frac{u_\lambda(t,x)}{v_\lambda(t,x)}=1
\eeq
for each $t\in\bbR$, at least for some $\lambda>\lambda_0$.  This has been achieved in \cite{NRRZ} for KPP reactions which are decaying (as $|x|\to\infty$) perturbations of a homogeneous KPP reaction, and for more general KPP reactions in \cite{ZlaInhomog}.  In both cases one needs $\lambda_0<2a_-$ (otherwise it is possible that no transition fronts exist \cite{MNRR}) and $\lambda\in(\lambda_0,2a_-)$.

In the present paper we show that this linearization approach can be extended to general {\it non-KPP monostable reactions.}  Our method is an extension of the (relatively simple and robust) approach from \cite{ZlaInhomog}.  There it was discovered that while $v_\lambda$ is obviously a super-solution of \eqref{1.1} when $g_1(u)=u$ (i.e., in the KPP case), one can also use $v_\lambda$ to find a sub-solution of the form $\til w_\lambda(t,x)=\til h_\lambda(v_\lambda(t,x))$, for $\lambda\in(\lambda_0,2a_-)$ and an appropriate $g_0$-dependent increasing function $\til h_\lambda: [0,\infty) \to [0,1)$ with 
\beq\lb{1.32}
\til h_\lambda(0)=0, \qquad \til h_\lambda'(0)=1,   \qquad \lim_{v\to\infty} \til h_\lambda(v)=1, \qquad  \til h_\lambda(v)\le v \quad \text{on $[0,\infty)$}.
\eeq
It follows that $\til w_\lambda\le v_\lambda$, and one can then find a transition front $u_\lambda$ between the two using parabolic regularity (see below).

While $\til w_\lambda$ remains a sub-solution for all $g_1$ from \eqref{1.3.2}, $v_\lambda$ need not be anymore a super-solution.  However, it turns out that one may still be able to construct a super-solution of the form $w_\lambda(t,x)=h_\lambda(v_\lambda(t,x))$, for an appropriate increasing $h_\lambda: [0,\infty) \to [0,\infty)$ such that
\beq\lb{1.33}
h_\lambda(0)=0, \qquad h_\lambda'(0)=1, \qquad  h_\lambda''(v)\ge 0 \quad \text{on $h_\lambda^{-1}([0,1])$}.
\eeq
Once again, we then find a transition front $u_\lambda$ between $\til w_\lambda$ and $\min\{w_\lambda,1\}$.

Moreover, a result of Nadin \cite{Nadin} (see also \cite{Shen}) shows that once some front exists, then also a (time-increasing) {\it critical front} exists.  The latter is a transition front $u_C$ for \eqref{1.1} such that
if $u\not\equiv u_C$ is any other transition front and $u(t,x)=u_C(t,x)$ for some $(t,x)\in\bbR^2$, then  
\[
[u_C(t,y)-u(t,y)](y-x)<0
\]
for all $y\neq x$.  That is, a critical front is the (unique up to time translation) ``steepest'' transition front for \eqref{1.1}, and is the inhomogeneous version of the minimal speed front for homogeneous reactions.  Indeed,  if $f$ is homogeneous, then $u_C$ is precisely the traveling front with the minimal speed $c_0$.

Thus we obtain the following result.

\begin{theorem}
\label{theorem:1.1}
Assume (\ref{1.30})--(\ref{1.4}), let $\nu:=\sup_{u\in(0,1]} \tfrac{g_1(u)}u  \ge 1$,
and let  the supremum of the spectrum of $\mathcal L:=\partial_{xx}+a(x)$ be $\lambda_0:=\sup \sigma(\mathcal L)\in[a_-,a_+]$.
If $\lambda\in(\lambda_0,2a_-)$ satisfies
\begin{equation} \lb{1.40}
 \lambda\le 2a_{-} - \frac{2\sqrt{\nu-1}}{\sqrt\nu +\sqrt{\nu-1}}a_{+},
\end{equation} 
then \eqref{1.1}  has a transition front $u_\lambda$ with $(u_\lambda)_t>0$, satisfying \eqref{1.41}.  In particular, if $\lambda_0$ is smaller than the right-hand side of \eqref{1.40}, then a  critical front $u_C$, with $(u_C)_t>0$, also exists.
\end{theorem}

{\it Remarks.}  1.  If $a(x)=f_u(x,0)$ is constant on $\bbR$, then $\lambda_0=a_-=a_+$, so the right-hand side of \eqref{1.40} is always greater than $\lambda_0$.  Thus a transition front exists for any $g_0,g_1$ in this case.
\smallskip


2.  The front $u_\lambda$ does not have a constant speed in general, but when $f$ is stationary ergodic in $x$, then it almost surely has an {\it asymptotic speed} $c_\lambda>0$ in the sense that if $X(t)$ is the rightmost point such that $u(t,X(t))=\tfrac 12$, then
\[
 \lim_{|t|\to\infty} \frac{X(t)}{t} =c_\lambda.
\]
\smallskip

3.  The result also holds with $v_\lambda$ replaced by more general solutions of \eqref{1.8} of the form $v_\mu(t,x)\equiv \int_\bbR  v_{\lambda}(t,x) d\mu(\lambda)$, with $\mu$ a finite non-negative non-zero Borel measure supported on a compact subset of $(\lambda_0,2a_- -  2\sqrt{\nu-1}(\sqrt\nu +\sqrt{\nu-1})^{-1}a_{+}]$ (or of $(\lambda_0,2a_-)$ if $\nu=1$).
\smallskip

4.  The result also applies to the more general equation 
\[
u_t = (A(x) u_x)_x + q(x) u_x + f(x,u)
\]
with 
\[
0<A_-\le A(x)\le A_+<\infty \qquad \text{and} \qquad  |q(x)|\le q_+<\infty
\]
for all $x\in\bbR$,  provided that $q_+\le 2\sqrt{(aA)_-}$ with $(aA)_-:=  \inf_{x\in\bbR} [a(x)A(x)]$, where 
\[
\lambda_0 := \sup_{\psi \in H^1(\bbR)} \frac{\int_\bbR [ - A(x)\psi'(x)^2 + q(x)\psi'(x)\psi(x) + a(x)\psi(x)^2] dx}{\int_\bbR \psi(x)^2 dx} \quad(\ge a_-)
\]
and $2a_-$ is replaced in \eqref{1.40} by
\[
\lambda_1:=\inf_{x\in\bbR} \left\{ a(x) +  \sqrt{(aA)_-} \left[ \sqrt{(aA)_-} - |q(x)| \right] A(x)^{-1} \right\} \quad(\le 2a_-).
\]
\smallskip

We indicate the proofs of Remarks 2--4 after the proof of the theorem.  

Our construction of the super-solution $w_\lambda$  is of independent interest and extends to more general equations in several dimensions, possibly with time-dependent coefficients.  Hence we state it here as a separate result.  

\begin{lemma} \lb{L.2.1}
Let the function $f(t,x,u)\ge 0$, positive definite matrix $A(t,x)$, and vector field $q(t,x)$ be all Lipschitz, with $(t,x,u)\in(t_0,t_1)\times \bbR^{d}\times[0,1]$ and some \hbox{$-\infty< t_0<t_1\le\infty$}.  Assume that $a(t,x)\equiv f_u(t,x,u)>0$ exists, \eqref{1.30}--\eqref{1.4} hold for all $(t,x,u)\in(t_0,t_1)\times\bbR^d\times[0,1]$, and define $\nu:=\sup_{u\in(0,1]} \tfrac{g_1(u)}u  \ge 1$.  Let $v>0$ be a solution of
\[
v_t = \nabla\cdot(A(t,x)\nabla v) + q(t,x)\cdot\nabla v + a(t,x) v
\]
on $(t_0,t_1)\times\bbR^d$.
If $\nu>1$ and for some $\alpha\le (\sqrt\nu-\sqrt{\nu-1})^2$ (or for some $\alpha<1$ if $\nu=1$),
\beq \lb{2.43}
\nabla v(t,x) \cdot A(t,x) \nabla v(t,x) \le \alpha a(t,x) v(t,x)^2
\eeq
holds for all $(t,x)\in (t_0,t_1)\times\bbR^d$, then there exist increasing functions $\til h$ satisfying \eqref{1.32} and $h$ satisfying \eqref{1.33} such that $\til w:= \til h(v)$ is a sub-solution  of 
\beq \lb{2.44}
u_t = \nabla\cdot(A(t,x)\nabla u) + q(t,x)\cdot\nabla u + f(t,x,u)
\eeq
on $(t_0,t_1)\times\bbR^d$ and $w:= h(v)$ is a super-solution on $[(t_0,t_1)\times\bbR^d]\cap\{(t,x)\,|\, w(t,x)\le 1\}$. 
Therefore, if $u$ solves \eqref{2.44} with 
\beq \lb{2.46}
\til w(t_0,x)\le u(t_0,x)\le \min\{w(t_0,x),1\}
\eeq
for all $x\in\bbR^d$, then for all $(t,x)\in (t_0,t_1)\times\bbR^d$ we have
\beq \lb{2.45}
\til w(t,x)\le u(t,x)\le \min\{w(t,x),1\}.
\eeq
\end{lemma}

\smallskip
{\bf Acknowledgements.}  All authors were supported in part   by the NSF grant DMS-1056327.  TT and BZ gratefully acknowledge the hospitality of the Department of Mathematics  at the University of Wisconsin--Madison during the  REU ``Analysis and Differential Equations'', where this research was performed.

\section{Proof of Theorem \ref{theorem:1.1} (using Lemma \ref{L.2.1})}

Let $\lambda\in(\lambda_0,2a_-)$ and $v=v_\lambda$ be from \eqref{1.31}, with $\phi=\phi_\lambda>0$ from \eqref{1.9} with $\lim_{x\to\infty} \phi(x)=0$ and $\phi(0)=1$.  It is proved in \cite{ZlaInhomog} that the unique such $\phi$  satisfies
\beq\lb{2.20}
\phi'(x)^2\le \alpha a(x)\phi(x)^2
\eeq
for $\alpha:=1-(2a_--\lambda)a_+^{-1}<1$ and all $x\in\bbR$, as well as
\beq\lb{2.21}
\phi(x)\ge 2\phi(y)
\eeq
for some $L<\infty$ and any $y-x\ge L$.

Since $\alpha\le (\sqrt\nu-\sqrt{\nu-1})^2$ is by the definition of $\alpha$ equivalent to
\[
\lambda\le 2a_--\left[1-\left(\sqrt\nu-\sqrt{\nu-1} \right)^2 \right] a_+ = 2a_{-} - \frac{2\sqrt{\nu-1}}{\sqrt\nu +\sqrt{\nu-1}}a_{+}
\]
(which is \eqref{1.40}),  Lemma \ref{L.2.1} applies to $v$ and \eqref{1.1}.  Thus we have \eqref{2.45} and a standard limiting argument now recovers an entire solution to \eqref{1.1} between $\til w$ and $\min\{w,1\}$.  We let $u_n$ be the solution of \eqref{1.1} on $(-n,\infty)\times\bbR$ with $u_n(-n,x):= \til w(-n,x)$.  Since $\til w(t,x)\le \min\{ w(t,x),1\}$ because $h(v)\ge v$ for $v\in h^{-1}([0,1])$, \eqref{2.46} is satisfied with $t_0:=-n$ 
and we have \eqref{2.45} on $(-n,\infty)\times\bbR$.
By parabolic regularity, there is a subsequence of $\{u_n\}$ which converges, locally uniformly on $\bbR^2$, to an entire solution $u$ of \eqref{1.1}.  We obviously have 
\beq\lb{2.25}
\til w\le u\le \min\{w,1\},
\eeq
and \eqref{1.41} for $u_\lambda:=u$ follows from $\til h'(0)=h'(0)=1$. 
We also have $u_t\ge 0$, since the same is true for $u_n$ due to $\til w_t=h'(v)v_t\ge 0$ and the maximum principle for $(u_n)_t$ (which satisfies a linear equation and is non-negative at $t=-n$).  The strong maximum principle then gives $u_t>0$ because obviously $u_t\not\equiv 0$. 
Finally, $u$ is a transition front because the second limit in \eqref{1.6} follows from $\lim_{x\to\infty} \phi(x)=0$ and \eqref{1.33}, and \eqref{1.7} holds with
\[
L_\eps:= L \left\lceil \log_2 \left( \til h^{-1}(1-\eps) -h^{-1}(\eps) \right) \right\rceil
\]
due to \eqref{2.25} and \eqref{2.21} (with $\til h,h$ from the lemma).  The first limit in \eqref{1.6} is then obvious from $u\le 1$, and the proof is finished by using the abovementioned result from \cite{Nadin} for critical fronts.
\smallskip

The claim in Remark 2 is proved as an analogous statement in \cite[Theorem 1.2]{ZlaInhomog}.

The claim in Remark 3 holds because  $L$ can be chosen uniformly for all $\lambda$ in the support of $\mu$ \cite{ZlaInhomog} and so \eqref{2.21} holds with $\phi(\cdot)$ replaced by $v_\mu(t,\cdot)$. Also, $v_\mu$ satisfies \eqref{2.20} with $\alpha$ corresponding to $\lambda:=\sup\supp\,\mu$.

The claim in Remark 4 holds because \eqref{2.20} and \eqref{2.21} continue to hold in that case, albeit with $2a_-$ replaced by $\lambda_1$ in the definition of $\alpha$ \cite{ZlaInhomog}.

\section{Proof of Lemma \ref{L.2.1}}

\cite{ZlaInhomog} shows that there is an increasing $\til h =\til h_\lambda$ as in \eqref{1.32} such that $\til w(t,x):=\til h(v)$ is a sub-solution of \eqref{1.1}.  This yields the first inequality in \eqref{2.45}.  We will next find an increasing $h=h_\lambda$ as in \eqref{1.33} such that  $w(t,x):=h(v(t,x))$ will be a super-solution where $w(t,x)\le 1$, which will yield the second inequality because $u\le 1$ by the hypotheses.  Our proof will be a super-solution counterpart to the sub-solution argument in \cite{ZlaInhomog}; it was a little surprising to us that such a counterpart argument exists for non-KPP reactions.

If $h$ is as in \eqref{1.33},  then \eqref{2.43} shows that
\begin{align*}
w_t -  \nabla\cdot(A\nabla w) - q\cdot\nabla w & = h'(v) [v_t -  \nabla\cdot(A\nabla v) - q\cdot\nabla v] - h''(v) \nabla v \cdot B \nabla v 
\\ & =h'(v) a v - h''(v) \nabla v \cdot B \nabla v 
\\ & \ge a [ vh'(v) - \alpha v^2 h''(v) ]
\end{align*}
when $w(t,x)\le 1$.  We can then conclude that $w$  is a super-solution of \eqref{1.1} where $w(t,x)\le 1$ once we show that on $h^{-1}([0,1])$ we also have
\beq \lb{2.22}
vh'(v)-\alpha v^2h''(v)\ge g_1(h(v)).
\eeq

It therefore remains to find $h$ satisfying \eqref{1.33} and \eqref{2.22}. 
We let $c:=\alpha^{1/2}+\alpha^{-1/2}$ and notice that since $\gamma+\gamma^{-1}\ge 2\sqrt\nu$ for all positive $\gamma\le \sqrt\nu-\sqrt{\nu-1}$, the hypothesis $\alpha\le (\sqrt\nu-\sqrt{\nu-1})^2$ yields $c\ge 2\sqrt\nu$.
Next let $U$ be the unique solution to the  ODE
\beq \label{2.3}
U''+cU' + g_1(U) = 0 
\eeq
on $[s_0,\infty)$, with
\beq \label{2.3a}
U(s_0)=1 \qquad\text{and}\qquad  U'(s_0)=-\sqrt\alpha g_1(1),
\eeq
where $s_0\in\bbR$ will be chosen later.  (This is the ODE that would be satisfied by the traveling front profile with speed $c$ for the homogeneous reaction $g_1(u)$ if we had $g_1(1)=0$; this profile would then also satisfy  $\lim_{s\to-\infty} U(s)=1$ and $\lim_{s\to-\infty} U'(s)=0$ instead of \eqref{2.3a}.) 

Notice that $U'(s_0)\ge -\tfrac c2$ because $g_1(1)\le\nu$ and 
\[
\sqrt\alpha \le  (\sqrt\nu+\sqrt{\nu-1})^{-1}\le \nu^{-1/2}.
\]
Let $V(s):=U'(s)$, and consider the curve $\tht:=\{(U(s),V(s))\}_{s\ge s_0}$.  It is easy to see that $\tht$ cannot leave the closed triangle $T$ in the $(U,V)$ plane with sides $V=0$, $U=1$, and $V=-\frac c2 U$.  This is because $(U(s_0),V(s_0))\in T$ and on $\partial T$, the vector field $(V,-cV-g_1(U))$ either points inside $T$ or is parallel to $\partial T$.  Here we use $c\ge 2\sqrt\nu$ to obtain on the third side
 \begin{equation} \label{3.3}
 \left(\frac{c}{2},1\right) \cdot \left( -\frac c2 U,-c\left(-\frac {c}2 U \right)-g_1(U) \right) = \frac{c^2}{4}U-g_1(U) \ge \nu U -g_1(U) \ge 0.
 \end{equation}
It follows that $U'(s)<0$ on $[s_0,\infty)$, and  since $g_1(U(s))> 0$, $U(s)$ cannot have local minima on $[s_0,\infty)$.
Hence  $\lim_{s\to\infty} U(s)$ exists and $\lim_{s\to\infty} U'(s)=0$. Finally, $g_1>0$ on $(0,1]$ yields $\lim_{s\to\infty} U(s)=0$.

%

We now define $h(0):=0$ and
\begin{equation} \label{2.1}
h(v) := U(-\alpha^{-1/2}\ln v)
\end{equation}
for $v\in(0,e^{-\sqrt\alpha s_0}]$, so $h$ is increasing and continuous at 0, with $h(e^{-\sqrt\alpha s_0})=1$ (we then extend $h$ onto $[0,\infty)$ arbitrarily, only requiring that it be increasing).
Since $c>2\sqrt{g_1'(0)}=2$ and
\[
\int_0^1 \frac{g_1(u)-u}{u^2}du <\infty
\]
by \eqref{1.3.1} and \eqref{1.4}, a result of Uchiyama \cite[Lemma 2.1]{Uchiyama} shows $\lim_{s \rightarrow \infty} {U(s)e^{\sqrt{\alpha}s}} \in (0,\infty)$.   (This result assumes $g_1(1)=0$ but we can extend $g_1, U$ to [0,2] so that $g_1(2)=0$ and $U$ satisfies \eqref{2.3}, and then apply \cite{Uchiyama} to $\til g(u):=\tfrac 12 g_1(2u)$ and the function $\til U(s):=\tfrac 12 U(s)$.) 

If we now pick the unique $s_0$ in \eqref{2.3a} such that $\lim_{s \rightarrow \infty} {U(s)e^{\sqrt{\alpha}s}} =1$ (notice that \eqref{2.3} is an autonomuous ODE), we obtain $h'(0)=1$.  We also have \eqref{2.22} on $[0,e^{-\sqrt\alpha s_0}]$ because on that interval,  \eqref{2.3} immediately yields 
\beq \lb{2.30}
\alpha v^2h''(v) - vh'(v)+ g_1(h(v))=0.
\eeq

It therefore  remains to show that  $h''(v)\ge 0$ on $[0,e^{-\sqrt\alpha s_0}]$.  Due to \eqref{2.30} and \eqref{2.1}, this is equivalent to
\beq\lb{2.31}
-U'(s) \ge \sqrt\alpha \, g_1(U(s))
\eeq
for $s\ge s_0$.  Thus we need to show that $\tht$ stays at or below  $\psi:=\{(U(s),-\sqrt\alpha\, g_1(U(s)))\}_{s\ge s_0}$.  This is true at $s=s_0$ by the definition of $U$, so it is sufficient to show that on $\psi$, the vector field $(V,-cV-g_1(U))$ points either below or parallel to $\psi$.  This holds because the normal vector to $\psi$ pointing down is $(\sqrt\alpha\,g_1'(U)V,V)$, so on $\psi$ we have
 \begin{eqnarray*}
(V,-cV-g_1(U))\cdot (\sqrt\alpha \,g_1'(U)V,V) 
 &=&\alpha^{1/2} g_1'(U)V^2-(\alpha^{1/2}+ \alpha^{-1/2})V^2- \alpha^{1/2} g_1(U)\alpha^{-1/2}V\\
 &=& \alpha^{1/2}(g_1'(U)-1)V^2,
 \end{eqnarray*}
which is non-negative due to \eqref{1.3.2} and $U(s)\le 1$ for $s\ge s_0$.  It follows that $h''(v)\ge 0$ on $[0,e^{-\sqrt\alpha s_0}]$ and so $h$ satisfies \eqref{1.33} and \eqref{2.22}.  The proof is finished.

\end{document}